\documentclass[12pt]{amsart}

\usepackage{amsfonts, amssymb, amscd}

\def\ZZ                 {{\mathbb Z}}
\def\QQ                 {{\mathbb Q}}
\def\PP                 {{\mathbb P}}
\def\RR                 {{\mathbb R}}
\def\CC                 {{\mathbb C}}
\def\nn                 {{\bf n}}
\def\pp                 {{\bf p}}
\def\vv                 {{\bf v}}

\newtheorem{lemma}{Lemma}[section]
\newtheorem{theorem}[lemma]{Theorem}
\newtheorem{corollary}[lemma]{Corollary} 
\newtheorem{proposition}[lemma]{Proposition}
\theoremstyle{definition}
\newtheorem{definition}[lemma]{Definition}

\newtheorem{remark}[lemma]{Remark}
\theoremstyle{remark}
\newtheorem*{proof*}{Proof}
\numberwithin{equation}{section}

\begin{document}
\title{Higher Stanley-Reisner rings and Toric Residues} 

\author{Lev A.  Borisov}
\address{Department of Mathematics \\ University of Wisconsin \\
  Madison \\ WI \\ 53706 \\ {\tt borisov@math.wisc.edu}.}

\thanks{The author was partially supported by NSF grant DMS-0140172.}

\begin{abstract}
We give a purely algebraic proof of the hypersurface case of Toric Residue
Mirror Conjecture recently proposed by Batyrev and Materov.
\end{abstract}

\maketitle
\section{Introduction}
Toric Residue Mirror Conjecture (TRMC) has been  formulated by 
Batyrev and Materov in \cite{BM}. It is in many ways analogous to
by now classical calculations of (virtual) numbers of rational curves
in Calabi-Yau hypersurfaces in toric varieties. In that story the 
generating function of the numbers of rational curves on a Calabi-Yau
hypersurface is calculated in terms of the periods of the mirror family.
Extensive references can be found in \cite{BM}.

Instead of  using Kontsevich's moduli spaces of stable curves on the 
ambient toric manifold in order to define virtual numbers of curves
on the hypersurface, the paper of Batyrev and Materov
uses a less sophisticated toric version of moduli 
spaces.  The resulting generating function is then conjectured to 
be related to the 
\emph{toric residue} of the mirror family. It is important to emphasize that
while the original mirror conjecture uses GKZ hypergeometric functions
of \cite{GKZhyper}, TRMC is formulated in terms
of some rational functions of several variables. This, perhaps, is the strongest
indication of the relative degree of difficulty of the two conjectures.
On the other hand, neither conjecture follows from the other.
Similar to the usual mirror symmetry Toric Residue Mirror 
Conjecture can be extended to the 
case of Calabi-Yau complete intersections defined by nef partitions, 
see \cite{BM2}.

In the present 
paper we give a simple algebraic proof of the hypersurface case of TRMC.
We build on the work done in \cite{BM} while at the same time
try to simplify it.
We do not attempt to use the most geometric version of the moduli
spaces of rational curves of given cohomology class on the toric variety,
but are willing to use a bigger space while adjusting the virtual fundamental
class on it. This allows us to essentially use a single cohomology space to
do all the calculations in. In fact, since 
we are mostly interested in cohomology classes, the moduli spaces we 
are working with do not have a direct geometric meaning.  However
we feel that they greatly simplify the exposition and allow for a
more conceptual understanding of TRMC. 

The paper is organized as follows. In Section \ref{sec.SR} we introduce
\emph{higher Stanley-Reisner rings}  $A_k$ of a toric variety, which is the main
tool of this paper. They are closely related to the cohomology of the 
moduli spaces considered
in \cite{BM} but are much easier to deal with. In Section \ref{sec.MP}
we mimic the construction of \cite{BM} to define Morrison-Plesser classes
in $A_k$ which are analogs of virtual fundamental classes on moduli
spaces of Kontsevich's stable curves. Section \ref{sec.hess} contains an 
explicit combinatorial calculation of the generating function of Morrison-Plesser
classes for a Hessian. This is the most delicate calculation of the entire paper.
Sections \ref{sec.SR}, \ref{sec.MP} and \ref{sec.hess} are self-contained
and can be read by anyone with just a minimum background in toric geometry.
By using somewhat technical results about secondary polytopes, 
we prove our version of TRMC in Section \ref{sec.TRMC}. Finally,
in Section \ref{sec.comp} we draw the connection between our definitions 
and that of \cite{BM}, thus establishing the hypersurface case of TRMC.

This is not the first solution of TRMC. In fact, the first proof of it belongs 
to Andr\'as Szenes and Mich\`ele Vergne. 
Although this paper is a result of a completely independent project,
the author has been informed by Victor Batyrev of the work of \cite{SV},
then still in preparation. Andr\'as Szenes then assured the author that
the two approaches differ sufficiently to warrant the completion of
the project. I thank both of them for their interest and encouragement.

\section{Higher Stanley-Reisner  rings}\label{sec.SR}
In this section we describe \emph{higher Stanley-Reisner rings} $A_k$ of
a toric variety $\PP_\Sigma$ for any positive integer $k$. The ring
$A_0$ is isomorphic to the cohomology ring of $\PP_\Sigma$, and
for any $k$ the ring $A_k$ admits a presentation inspired by Stanley-Reisner
description of $A_0$. This is the reason behind our terminology.

Let $\Sigma$ be a complete simplicial fan in a lattice $M\cong \ZZ^d$ and let 
$v_i,i=1,\ldots,n$ be the minimum generators of its one-dimensional cones.
Let $\PP_\Sigma$ be the corresponding toric variety. Its cohomology
is given by the Stanley-Reisner relations:
$$
H^*(\PP_\Sigma,\CC)\cong\CC[D_1,\ldots,D_n]/I
$$
where the ideal $I$ is generated by linear relations 
$\sum_{i=1}^n (\lambda\cdot v_i) D_i$ for all $\lambda\in N=M^*$ and monomial relations
$\prod_{i=1}^n D_i^{r_i}$
over all $\{r_i\}$ such that \emph{no cone of $\Sigma$ contains all 
$v_i$ for which $r_i>0$}. This is a slight reformulation of 
the usual description, which the reader can easily verify to be
equivalent.

\begin{definition}
For every nonnegative $k$ we denote by $A_k$ the quotient of the polynomial
ring $\CC[D_1,\ldots,D_n]$ 
by linear relations $\sum_{i=1}^n (\lambda\cdot v_i) D_i $, $\lambda\in N$
and monomial relations $\prod_{i=1}^nD_i^{r_i}$ 
over all $\{r_i\}$ such that \emph{no cone of $\Sigma$ contains all 
$v_i$ for which $r_i>k$}.
\end{definition}

We will now show that $A_k$ is isomorphic to the cohomology of 
some complete toric variety of dimension $nk+d$, defined by a fan
$\Sigma_k$ in the lattice $M\oplus\ZZ^{nk}$, as follows.
Let $\{e_{i,j},1\leq i\leq n, 1\leq j\leq k\}$ be a basis in $\ZZ^{nk}$.
For each $i$ we introduce $e_{i,0}=-\sum_{j=1}^k e_{i,j}$.
We then consider elements 
$$v_{i,j} = v_i \oplus e_{i,j},~1\leq i\leq n,0\leq j\leq k.
$$ 
Cones of the fan $\Sigma_k$ are generated by collections
of elements $v_{i,j}$ such that the indices $i$ that occur $(k+1)$ times 
correspond to generators $v_i$ of a cone of $\Sigma$. 

\begin{proposition}\label{fan}
The above described $\Sigma_k$ is a complete simplicial fan. 
The cohomology ring of the corresponding toric variety 
$\PP_{\Sigma_k}$ is isomorphic
to $A_k$.
\end{proposition}

\begin{proof}
First of all, we need to see that $\Sigma_k$ is a fan, i.e. the intersection
of two cones $C_1$ and $C_2$ in it is again a cone in $\Sigma_k$.
It is sufficient to show that if $C_i$ correspond to the subsets of 
indices $I$ and $J$ of $\{1,\ldots,n\}\times\{0,\ldots,k\}$ then 
$C_1\cap C_2$ is equal to the cone $C$ spanned by $v_{i,j}$ for 
$(i,j)\in I \cap J$.  For each $i=1,\ldots, n$ we denote by $I_i$, 
$J_i$ and $(I\cap J)_i$ the $i$-th components of $I$, $J$ 
and $I\cap J$ respectively.

It is clear that $C\subseteq C_1\cap C_2$.
To show the converse suppose that $w=v\oplus \bigoplus_{i=1}^nw_i$
is in $C_1\cap C_2$. 
We have 
$$
w= \sum_{i=1}^n (\sum_{j\in I_i} \alpha_{i,j}) v_i
\oplus \bigoplus_{i=1}^n \sum_{j \in I_i} \alpha_{i,j}e_{i,j}
=\sum_{i=1}^n (\sum_{j\in J_i} \beta_{i,j}) v_i
\oplus \bigoplus_{i=1}^n \sum_{j \in J_i} \beta_{i,j}e_{i,j}
$$
where all $\alpha$ and $\beta$ are nonnegative.
For each $i$ we have $\alpha_{i,j} = \beta_{i,j} + \gamma_i$ 
for some numbers $\gamma_i$ independent of $j$. We observe 
that if $\gamma_i>0$ then $\alpha_{i,j}>0$ for all $j$, so
$\vert I_i \vert=k+1$. Similarly, if 
$\gamma_i<0$ then $\vert J_i\vert =k+1$. We have
$$
0=\sum_{i=1}^n  \gamma_i v_i.
$$
By splitting this into the sums with positive and negative $\gamma_i$
we get 
$$
\sum_{i, \vert I_i \vert = k+1,\gamma_i>0} \gamma_i v_i = 
\sum_{i, \vert J_i \vert = k+1,\gamma_i<0}  (-\gamma_i) v_i.
$$
By the definition of cones in $\Sigma_k$, the set of $v_i$ with 
$\vert I_i \vert = k+1$ forms a cone in $\Sigma$, and similarly for $J$.
Both sides of the above identity lie in the intersection of the 
corresponding cones. As a result, $\gamma_i=0$ unless 
$\vert I_i \vert = \vert J_i \vert = k+1$. This implies that $\alpha_{i,j}
=\beta_{i,j}$ unless $\vert I_i \vert = \vert J_i \vert = k+1$. Consequently,
if $(i,j)\not\in J$ then $\vert J_i\vert <k+1$ so $\alpha_{i,j}=\beta_{i,j}=0$.
Hence only nonzero $\alpha_{i,j}$ come from $(i,j)\in I\cap J$.
This shows that $\Sigma_k$ is indeed a fan.

To show that $\Sigma$ is complete, consider any
$w=v\oplus \bigoplus_{i=1}^n w_i$. Each $w_i$ sits in the unique cone 
of the standard fan for $\PP^k$ so it can be written in a unique way as
a nonnegative linear combination of $k$ vectors $e_{i,j}$. If we subtract
the corresponding linear combinations of $v_{i,j}$ we are left with 
an element $v'$ of $M_\RR$. We can write it as a positive linear combination
 $v'=\sum_{v_i\in \sigma} \gamma_i v_i$ for some cone $\sigma\in \Sigma$.
Therefore, $v'\oplus{\bf 0}=\sum_{v_i\in\sigma} \sum_{j=0}^k
\frac 1{k+1} \gamma_i v_{i,j}$. It is easy to see that the resulting linear
combination for $w$ will have positive coefficients for the set of indices $I$ such that $\vert I_i\vert= k+1$ iff $v_i\in \sigma$. Thus $w$ is lies in a cone of $\Sigma_k$.

To calculate the cohomology of $\PP_{\Sigma_k}$ we use the Stanley-Reisner
presentation of it as a quotient of a polynomial ring in $n(k+1)$ variables
$D_{i,j}$ by linear and polynomial relations. We have  linear relations 
$D_{i,j_1}=D_{i,j_2}$ for all $i$, $j_1$ and $j_2$, which come from
the linear functions on each copy of $\ZZ^k$. 

We can  map $D_i$ to $(k+1)D_{i,0}$  and note that linear relations 
coming from $M$ give
$$
\sum_{i=1}^n (\lambda\cdot v_i) D_i = 0.
$$
The description of the cones of $\Sigma_k$ then shows that the monomial 
relations are exactly the ones in the definition of $A_k$.
\end{proof}

\begin{remark}\label{shift}
For any  $l>0$ the ring $A_k$ can be naturally mapped to
$A_{k+l}$ by multiplying by $\prod_{i=1}^n D_i^l$. Indeed, multiplication by $\prod_{i=1}^n D_i^l$ maps monomial relations for $A_k$ into monomial relations for $A_{k+l}$ so it maps the ideal of relations for $A_k$ into
that for $A_{k+l}$.  This map is a $\CC[D_1,\ldots,D_n]$-module 
map. In what follows we will be  working in the direct limit
of $A_k$ under these maps. 
\end{remark}

By Proposition \ref{fan}, the component of degree $nk+d$ in the 
graded ring $A_k$ is one-dimensional. Moreover we have evaluation
maps $\int_{\PP_{\Sigma_k}}\colon A_k\to\CC$ coming the intersection 
on $\PP_{\Sigma_k}$. We observe that these evaluations are not quite compatible
with the maps of Remark \ref{shift}.

\begin{proposition}\label{s}
For any element $a\in A_k$ and any $l>0$ there holds
$$
\int_{\PP_{\Sigma_k}}a = 
\frac
{(k+1)^{nk}}
{(k+l+1)^{n(k+l+1)}}
\int_{\PP_{\Sigma_{k+l}}} a \, \prod_{i=1}^n D_i^l.
$$
\end{proposition}

\begin{proof}
Since the top degree components of $A_k$ are one-dimensional,
the above statement is true up to a multiplication by a constant.
Consequently, it is enough to show that 
\begin{equation}\label{int}
\int_{\PP_\Sigma} a = (k+1)^{-nk} \int_{\PP_{\Sigma_k}} 
a \, \prod_{i=1}^n D_i^k
\end{equation}
for one nonzero element $a$ of $A_0$. Pick a maximum
cone $\sigma\in\Sigma$ and let $V(\sigma)$ be the \emph{normalized
volume} of the simplex generated by $v_i\in \sigma$, which is defined
as the absolute value of the determinant of $v_i\in \sigma$ expanded 
in a basis of the lattice. Then 
$$
\int_{\PP_\Sigma}\prod_{v_i\in \sigma}D_i= V(\sigma)^{-1}.
$$
On the other hand, 
$\prod_{v_i\in\sigma} D_i \prod_{i=1}^n D_i^k$
corresponds to 
$$
\prod_{v_i\in\sigma}(k+1)D_{i,0}\,\prod_{i=1}^n\prod_{j=1}^k(k+1)D_{i,j}.
$$
We need to find the normalized volume of the cone
$\sigma_k$ in $\Sigma_k$ generated by 
$v_{i,0}$ for $v_i\in\sigma$ and $v_{i,j}$ by $1\leq i\leq n, 1\leq j\leq k$
in $M\oplus \ZZ^{nk}$. We claim it equals $V(\sigma)(k+1)^d$. Indeed, 
replacing $v_{i,0}$ by $\sum_{j=0}^kv_{i,j}=(k+1)v_i\oplus {\bf 0}$ does not 
change the volume. Then it is easy to calculate the volume of the resulting 
cone. Then we have 
$$
\int_{\PP_{\Sigma_k}} 
\prod_{v_i\in \sigma}D_i\, \prod_{i=1}^n D_i^k=(k+1)^{nk+d}V(\sigma_k)^{-1}
=(k+1)^{nk}V(\sigma)^{-1}.
$$
This shows \eqref{int}.
\end{proof}

We will adjust the top class evaluation so that it is compatible with the 
maps of Remark \ref{shift}.
\begin{definition}
For each nonnegative $k$ we define $\int_{A_k}\colon A_k\to \CC$ by
$\int_{A_k} = (k+1)^{-nk} \int_{\PP_{\Sigma_k}}$. By Proposition \ref{s},
we have 
\begin{equation}\label{sameint}
\int_{A_k} a = \int_{A_{k+l}} 
a \, \prod_{i=1}^n D_i^l
\end{equation}
for all $k,l$ and all $a\in A_k$.
\end{definition}

\begin{definition}\label{intma}
We denote by $\mathcal A$ the direct limit of $A_k$ taken with respect
to the maps of Remark \ref{shift}. The direct limit of $\int_{A_k}$ gives 
a map $\int_{\mathcal A}
\colon{\mathcal A}\to \CC$. In addition $\mathcal A$ inherits
a structure of $\CC[D_1,\ldots,D_n]$-module.
\end{definition}

\begin{remark}\label{nondiag}
Arguments of this section do not depend on the fact that the same lattice $\ZZ^k$
is used for all $i$. In fact, for any nonnegative integers $k_1,\ldots,k_n$ one
can define the fan $\Sigma_{(k_1,\ldots,k_n)}$ in the lattice
$$
M\oplus\bigoplus_i \ZZ^{k_i},
$$
in terms of  $v_{i,j}=v_i\oplus e_{i,j}$ for $1\leq i\leq n,0\leq j\leq k_i$.
The cohomology ring $A_{(k_1,\ldots,k_n)}$
of the corresponding toric variety $\PP_{\Sigma_{(k_1,\ldots,k_n)}}$ 
is given by the usual linear relations and the relations $\prod_i D_i^{r_i}=0$ 
if the set of $v_i$ for which $r_i>k_i$ does not lie in a cone of $\Sigma$. 
Products of powers of $D_i$ define maps between these rings, which 
are compatible with the top class evaluation, once it is adjusted by 
the factor $\prod_{i=1}^n (k_i+1)^{k_i}$. While introducing these rings
has no effect on the limit $\mathcal A$, we will use this remark in
Section \ref{sec.comp}.
\end{remark}

The following convention will be used in the later sections.
\begin{definition}
We define $D_0=-\sum_{i=1}^n D_i$ in each ring $A_k$.
\end{definition}

\section{Morrison-Plesser classes}\label{sec.MP}
In what follows it will be convenient to extend the lattice $M$ to
a lattice $\bar M: = M\oplus \ZZ$. We will denote by $\vv_k$
the elements $v_k\oplus 1$ of $\bar M$. We will also consider
$\vv_0={\bf 0}\oplus 1\in \bar M$. In what follows we will consider
linear combinations of $\vv_k$, which are encoded by elements
$\beta=(b_0,b_1,\ldots,b_n)\in \ZZ^{n+1}$. From this section on
we assume that the toric variety $\PP_\Sigma$ is nef-Fano.
This means that $v_i$ lie on the boundary of the convex polytope $\Delta
={\rm conv}(\{v_i,i=1,\ldots,n\})$.

\begin{definition}
Let $\beta\in \ZZ^{n+1}$ be any lattice point which satisfies
$b_0\leq 0$. 
We define the \emph{Morrison-Plesser class}
as an element of $\mathcal A$ which is the image of 
the element in $A_k$
$$
\Phi_{\beta} = (D_0)^{-b_0} \prod_{i=1}^n D_i^{k-b_i}.
$$
for some sufficiently big $k$. Clearly, the result is independent of
a choice of $k$.
\end{definition}

The following key construction is motivated by \cite{BM}.
Let $f=1+\sum_{i=1}^n a_it^{v_i}$ be a generic formal Laurent polynomial
in $t^M$. Notice the change of sign in our notations as compared to 
that of \cite{BM}.
We will also use notation $a_0=1$.

Denote by $K$ the cone in $\bar M_\RR$ spanned by $\vv_k$. It can be also
described as $\{c\Delta\oplus c, c\geq 0\}$.
\begin{definition}
Let $\pp\in \bar M$
be a point in $K$. Then we define a formal Laurent series 
$$
\Psi_\pp :=
\sum_{\beta:\sum_{i}b_i\vv_i=-\pp, b_0\leq 0} \Phi_{\beta}\prod_{i=1}^n
a_i^{b_i}
$$
with values in $\mathcal A$.
\end{definition}

\begin{remark}
The above definition gives $\Psi_\pp=0$ if $\pp$ does not lie in 
the lattice generated by $\vv_i,i=0,\ldots,n$. 
\end{remark} 

\begin{proposition}\label{module}
For every $\pp\in K$ and every $j\in\{0,\ldots,n\}$ there holds
$$
D_j\Psi_\pp = \alpha_j\Psi_{\pp+\vv_j}.
$$
\end{proposition}

\begin{proof}
The equality has to be understood as that of formal Laurent series
in $a_1,\ldots,a_n$. 

Every solution of $\sum_{i}b_i\vv_i=-\pp$ gives a solution
$\sum_i\hat b_i\vv_i=\sum_i (b_i-\delta_i^j)\vv_i=-\pp-\vv_j$ and vice versa. 
For $j>0$ the coefficients $b_0$ in these two solutions are the same,
and it is straightforward to see that the elements of $\mathcal A$ for 
the left and right hand sides of the above equation are the same.

The situation is more complicated in the case of $j=0$. The 
summation for $\Psi_{\pp+\vv_0}$ involves solutions with
$\hat b_0 =0$ which have no counterpart in the summation for
$\Psi_\pp$. We will however see that these elements are in fact
zero in $\mathcal A$. We will do the calculation in $A_k$ for all $k$ 
for which a given $\Phi_{\hat \beta}$ makes sense.

Suppose that  $\Phi_{\hat\beta}$ is not zero. It is proportional
to the monomial 
$$
\prod_{i=1}^n D_i^{k-\hat b_i},
$$
so the definition of $A_k$ implies that the set of $v_i$ for which 
$\hat b_i<0$ lies in some cone $\sigma\in\Sigma$. Consider the corresponding
face in $K$, generated by all $\vv_i$ for which $\hat b_i<0$. Consider
its supporting hyperplane given by an element $\nn$ of $\bar N$. 
We have $\nn\cdot \vv_i\geq 0$ for all $1\leq i\leq n$ and
$\nn\cdot \vv_i =0$  if $\hat b_i<0$. Since we have 
$$\sum_{i=1}^n \hat b_i\vv_i = -\pp-\vv_0$$
we conclude $\nn\cdot (\pp+\vv_0)\leq 0$. However, by assumption
$\pp\in K$, so $\pp+\vv_0$ is in the interior of $K$ and every supporting
hyperplane is strictly positive on it.
\end{proof}

\begin{proposition}\label{rels}
For every $\pp\in K$ and every element $\nn\in \bar N=\bar M^*$ there holds
$$
\sum_{i=0}^n a_i(\nn\cdot \vv_i) \Psi_{\pp+\vv_i} = 0.
$$
\end{proposition}

\begin{proof}
By Proposition \ref{module} we get 
$$
\sum_{i=0}^n a_i(\nn\cdot \vv_i) \Psi_{\pp+\vv_i} 
=\sum_{i=0}^n (\nn\cdot\vv_i)D_i\Psi_{\pp}.
$$
For $\nn=\lambda\oplus 0$ the result follows from the linear relations
that hold in each $A_k$, hence in $\mathcal A$. 
For $\nn={\bf 0}\oplus 1$ the
result follows from the definition of $D_0$. By linearity, the statement
holds for all $\nn$.
\end{proof}

We will be especially interested in $\Psi_\pp$ of top degree.
\begin{definition}
Consider points $\pp\in K$ given by $p\oplus d\in \bar M\cong M\oplus \ZZ$. For each such $\pp$ we define a formal Laurent series in
$a_i$ by
$$
\int_{\mathcal A}\Psi_p:=
\sum_{\beta:\sum_{i}b_i\vv_i=-\pp, b_0\leq 0} 
\int_{\mathcal A}\Phi_{\beta}\prod_{i=1}^n
a_i^{b_i}
$$
where $\int_{A_k}$ is the top class evaluation in 
$\mathcal A$, see Definition \ref{intma}. 
\end{definition}

\section{Hessians}\label{sec.hess}
The goal of this section is to calculate the linear combinations of series 
$\int_{\mathcal A}\Psi_\pp$ for linear combinations of points $\pp\in K$ that come from certain Hessians. We will now assume that we have 
a reflexive polytope $\Delta\subset M$ and let $\mathcal T$ be
a triangulation of $\Delta$ whose maximum simplices contain $\bf 0$.
The vertices of $\mathcal T$ are $\{{\bf 0}\} \cup\{v_i,i=1,\ldots,n\}$. Points $v_i$
include all vertices of $\Delta$ and
perhaps some other points $v_i\in\partial\Delta$. 
The triangulation $\mathcal T$ induces a complete simplicial fan $\Sigma$ in $M_\RR$. As before, we introduce $\bar M\cong M\oplus\ZZ$, $\vv_i$
and $K$. The reflexivity assumption on $\Delta$ implies that 
every lattice point in the interior of $K$ lies in $\vv_0 + K$.

As before we fix a generic polynomial $f=a_0 +\sum_{i=1}^na_iv_i$
with $a_0=1$. 
We recall the definitions of  Hessians 
$H_f$ and $H'_f$ associated to these data, see \cite{BM}.
We will give a mostly self-contained exposition, both for the benefit of 
the reader and to facilitate further arguments.

\begin{definition}
Fix a basis $\{\nn_1,\ldots,\nn_{d+1}\}$ of the lattice $\bar N=\bar M^*$.
The Hessian $H_f$ is the determinant of the square matrix of size $d+1$
whose $(i,j)$-th entry is 
$$
\sum_{l=0}^n (\nn_i\cdot \vv_l)(\nn_j\cdot \vv_l) a_l t^{\vv_l}
$$
where $t$ is a dummy variable. 
\end{definition}

\begin{remark}\label{detsq}
The Hessian $H_f$ can be thought of as an element of $\CC[K]$.
It is easy to see that $H_f$ does not depend on the choice of the
basis of $\bar N$. Indeed, any linear change in $\nn_i$ by a matrix $R$
amounts to the linear change on rows and columns of the above matrix
and gives an additional factor of $(\det R)^2$ to the Hessian.
\end{remark}

\begin{proposition}\label{formula}\cite{CDS}
There holds
$$
H_f = \sum_{J\subseteq \{0,\ldots,n\},
\vert  J\vert =d+1} V(J)^2
\,
\big(\prod_{i\in J} a_i\big)
\,
t^{\sum_{i\in J}\vv_i} 
$$
where $V(J)$ is the normalized $(d+1)$-dimensional volume 
of the simplex spanned by the vectors $\vv_i,$ for $i\in J$.
\end{proposition}

\begin{proof}
It is clear that all terms of the determinant involve at most 
$(d+1)$ different $a_it^{\vv_i}$. Consequently, it is enough to find
out what happens when all $a_i$ are zero except for $i\in J,\vert J\vert =d+1$. 

If the elements $\vv_i,i\in J$ are linearly dependent, there is an element of
$\nn_1\in \bar N$ that vanishes on all of them. By completing it to
the basis we see that the Hessian is the determinant of a matrix
with a zero first row (and zero first column). Consequently, these 
collections $J$ do not contribute to the Hessian. 

If the elements $\vv_i,i\in J$ are linearly independent, consider the dual
basis of $\bar N_\QQ$. The matrix will then be simply the diagonal
matrix with $(i,i)$-th entry $a_it^{\vv_i}$. The dual basis will typically
not be a basis of $\bar N$. However, a basis of $\bar N$ is obtained
by a linear transformation of determinant $V(J)$ from the dual basis 
to $\vv_i,i\in J$. The argument of Remark \ref{detsq} then completes the proof.
\end{proof}

As a consequence of the above proposition, $H_f$ is supported in
the interior of $K$. Indeed, the terms from the boundary correspond
to $J$ with $V(J)=0$. Since $\Delta$ is reflexive, $H_f$ is divisible
by $t^{\vv_0}$ in $\CC[K]$ which allows us to introduce 
$H'_f$.
\begin{definition}
$$
H'_f:=H_f/t^{\vv_0}=\sum_{J\subseteq \{0,\ldots,n\},
\vert  J\vert =d+1} V(J)^2
\,
\big(\prod_{i\in J} a_i\big)
\,
t^{\sum_{i\in J}\vv_i-\vv_0} 
$$
\end{definition}

The main result of this section is the following calculation which describes the
value of $\int_{\mathcal A}\Psi$ on the Hessian. We denote by 
${\rm Vol}(\Delta)$ the normalized volume of $\Delta$.
\begin{theorem}\label{hess}
$$
\sum_{J\subseteq \{0,\ldots,n\},
\vert  J\vert =d+1} V(J)^2
\,
\big(\prod_{i\in J} a_i\big)
\,
\int_{\mathcal A}\Psi_{\sum_{i\in J}\vv_i-\vv_0}  = {\rm Vol}(\Delta).
$$
\end{theorem}

\begin{proof}
The definition of $\Psi_{\sum_{j\in J}\vv_j-\vv_0}$ involves the summation 
over $\hat \beta=(\hat b_0,\ldots,\hat b_n)$ with $\hat b_0\leq 0$ and
$$\sum_{i=0}^n \hat b_i \vv_i  = -\sum_{j\in J}\vv_j + \vv_0.$$
We introduce $b_i=\hat b_i + \chi(i\in J) - \delta_i^0$. Here 
$\chi(i\in J)$ is $1$ if $i\in J$ and is zero otherwise and $\delta$ is the
Kronecker symbol.
Then we have
the sum over $\beta = (b_0,\ldots,b_n)$ with 
$$
\sum_{i=0}^n b_i\vv_i = {\bf 0}
$$
and the additional assumption $b_0-\chi(0\in J)+1\leq 0$. This means
that the sum takes place over all $\beta$ with $b_0\leq  0$, but
for $b_0=0$ one only uses $J$ that contain $0$. 

We will analyze the contribution of $\beta$ from the following three
cases: $b_0<0$; $b_0=0,\beta\neq 0$; $\beta=0$. We will establish
the claim of the Proposition by showing that the only nonzero contribution
comes from $\beta=0$ and equals ${\rm Vol}(\Delta)$.

\emph{Case $b_0<0$.} The contribution is given by 
$$
\int_{A_k} \sum_{J\subseteq \{0,\ldots,n\},\vert J\vert =d+1} V(J)^2
\prod_{i\in J} a_i
\prod_{i\in J}  
D_0^{-\hat b_0} \prod_{i=1}^n \big( D_i^{k-\hat b_i}a_i^{\hat b_i}\big)
$$
$$
=\int_{A_k} D_0^{- b_0-1} \prod_{i=1}^n \big( D_i^{k-b_i}a_i^{b_i}\big)
\sum_{J\subseteq \{0,\ldots,n\},\vert J\vert =d+1}
V(J)^2 \prod_{i\in J} D_i.
$$
The proof of Proposition \ref{formula} shows that 
$\sum_{J\subseteq \{0,\ldots,n\},\vert J\vert =d+1}
V(J)^2 \prod_{i\in J} D_i$ is the determinant of square matrix of size $(d+1)$
whose $(i,j)$-th entry is
$$
\sum_{l=0}^n (\nn_i\cdot \vv_l)(\nn_j\cdot \vv_l) D_l.
$$
Here $\{\nn_i\}$ is an arbitrary basis of $\bar N$ so we can pick $\nn_1$
to have $\nn\cdot\vv_l=1$ for all $l$. Then the first row of the matrix
consists of elements that are zero in $A_k$, so the determinant is zero.
As a consequence, elements with $b_0<0$ do not contribute to the 
overall sum.

\emph{Case $b_0=0,\beta\neq 0$}.
This is the most difficult part of the calculation. We again would like to 
show that the contribution is zero. We recall that we have a summation
over the subsets $J$ that contain $0$. We abuse notations and use the
same letter for the corresponding subset of $\{1,\ldots,n\}$.
We need to show that 
\begin{equation}\label{need}
\int_{A_k} \prod_{i=1}^n D_i^{k-b_i}
\sum_{J\subseteq \{1,\ldots,n\},\vert J\vert =d}
V(J)^2 \prod_{i\in J} D_i=0
\end{equation}
where $V(J)$ now denotes the normalized volume of the $d$ vectors
$v_i,i\in J$ in the lattice $M$.

We immediately observe that \eqref{need} holds unless the set of
$v_i$ such that $b_i<0$ lies in a cone in $\Sigma$. Indeed, otherwise
we have $\prod_{i=1}^n D_i^{k-b_i}=0$ in $A_k$. We will denote the 
cone spanned by $v_i$ with $b_i<0$  by $\sigma$. We denote by
$\theta$ the minimum face of the reflexive polytope $\Delta$ that
contains all $v_i$ with $b_i<0$. Since we have 
$$
\sum_{i=1}^n b_i\vv_i = {\bf 0},
$$
all nonzero $b_i$ correspond to elements $v_i\in \theta$. Indeed,
there is an element $\nn$ in $\bar N$ which vanishes on $\vv_i$
for $v_i\in \theta$ and is positive on all other $\vv_i$. When applied
to both sides of the above equation we see that $\sum_{v_i\not\in\theta}
b_i (\nn\cdot \vv_i) = 0$. This shows that all terms $b_i (\nn\cdot \vv_i) $
are zero, since
all terms are nonnegative.

The second observation
is that $\prod_{i=1}^n D_i^{k-b_i}\prod_{i\in J} D_i$ is zero in $A_k$
unless all $v_i$ for $i\in J$ lie in a codimension one face 
$\theta_1\subset \Delta$ that contains $\theta$. Indeed, the set of 
exponents that 
are bigger than $k$ contains all $v_i\in \sigma$ and all
$v_i\not\in \theta, i\in J$. To be nonzero in $A_k$ implies that all these
elements lie in a cone of $\Sigma$. Hence the minimum face $\theta_2$ in 
$\Delta$ that contains these elements has codimension at least one.
Since this minimum face contains $\theta$, in fact all elements 
$v_i,i\in J$ lie in $\theta_2$. Then any codimension one face $\theta_1
\supseteq \theta_2$ works. In fact,  for nonzero $V(J)$ the face 
$\theta_1$ is uniquely determined. As a result, we can split the summation
over all $J$ into sub-summations over $\theta_1$. 
Then \eqref{need} would follow from
\begin{equation}\label{need2}
\int_{A_k} \prod_{i=1}^n D_i^{k-b_i}
\sum_{J\subseteq vert(\theta_1),\vert J\vert =d}
V(J)^2 \prod_{i\in J} D_i=0
\end{equation}
where $vert(\theta_1)$ is the set of indices $i$ for which $v_i\in \theta_1$.

For any basis $\{\lambda_1,\ldots,\lambda_d\}$ of $N$ consider a
square matrix $B$ with entries
$$
B_{ij}=\sum_{v_l\in \theta_1} (\lambda_i\cdot v_l)(\lambda_j\cdot v_l) D_l.
$$
Similarly to the proof of Proposition \ref{formula}, we can see
that 
$$
\sum_{J\subseteq vert(\theta_1),\vert J\vert =d}
V(J)^2 \prod_{i\in J} D_i={\rm Det} (B).
$$
Since we are only trying to show that this determinant is zero,
we could use a basis of $N_\QQ$ instead of $N$.
We will pick a special basis as follows. Element $\lambda_1$ will
be equal to $1$ on all $v_i\in \theta_1$. Elements 
$\lambda_{2},\ldots,\lambda_{r}$ 
where $r=\dim(\theta_1)-\dim(\theta)+1$
will be zero on $\theta\subseteq \theta_1$. It is easy to see that these
elements could be completed to a basis, if $\theta$ is non-empty
which is guaranteed by $\beta\neq 0$.

The first row of the matrix $B$ consists of $B_{1,j}=
\sum_{v_l\in \theta_1} (\lambda_j\cdot v_l) D_l$. These elements 
equal $-\sum_{v_l\not\in\theta_1} (\lambda_j\cdot v_l) D_l$ due
to linear relations in $A_k$. We are going to 
replace the first row of $B$ by the above elements of $\CC[D]$
and call the resulting matrix $B'$. We then calculate the determinant
of $B'$ as an element in $\CC[D]$.
We claim that all monomials in $D_i$ that appears in the resulting 
expression do not have $v_i$ lie in \emph{any} face of $\Delta$ that 
contains $\theta$. Consequently, the above arguments show that their contribution to the left hand side of \eqref{need2} are zero.

To substantiate our claim, we expand the determinant 
of $B'$ along the first $r$ rows. It is sufficient to show that 
all the $r\times r$ minors of the first $r$ rows of $B'$ have nonzero
coefficients only by monomials 
$$
\prod_{i=1}^r D_{l_i}
$$
such that no proper face $\theta_2\supseteq \theta$ contains
all of $v_{l_i}$.  
Suppose such monomial and such $\theta_2$ exist.
As in the proof of Proposition \ref{formula} we can replace
the rows $2,\ldots, r$ of $B'$ by keeping only the linear combinations
of $D_{l_i}$. We call the resulting matrix $B''$.
The face $\theta_2$ can not equal 
$\theta_1$, since the first row of $B'$ has $D_l$ with $v_l\not\in\theta_1$.
The intersection of $\theta_1$ and $\theta_2$ is a proper subface of 
$\theta_1$.  Consequently, there is a linear combination of $\lambda_2,
\ldots,\lambda_r$ which vanishes on $\theta_1\cap\theta_2$. By taking 
the appropriate linear combination of the rows $2,\ldots,r$ of $B''$,
we get a zero row, which means that the monomial 
$\prod_{i=1}^r D_{l_i}$ occurs with zero coefficient.

\emph{Case $\beta=0$.} This case is essentially covered in \cite{BM}
but we reproduce the argument here. The contribution equals
$$
\int_{A_k} \prod_{i=1}^n D_i^{k}
\sum_{J\subseteq \{1,\ldots,n\},\vert J\vert =d}
V(J)^2 \prod_{i\in J} D_i=\int_{\PP_\Sigma} \prod_{i=1}^n 
\sum_{J\subseteq \{1,\ldots,n\},\vert J\vert =d}
V(J)^2 \prod_{i\in J} D_i
$$$$=\sum_{\sigma\in \Sigma,\dim\sigma=d}
V(\sigma)^2\int_{\PP_\Sigma}\prod_{v_i\in \sigma} D_i =
\sum_{\sigma\in \Sigma,\dim\sigma=d}
V(\sigma)={\rm Vol}(\Delta).
$$
Here $V(\sigma)$ is the normalized volume of the corresponding simplex
of the triangulation.
This finishes the proof of Theorem \ref{hess}.
\end{proof}

\section{Toric Residue Mirror Conjecture}\label{sec.TRMC}
As before, we are working with a reflexive polytope $\Delta$,
a subset $\{v_i\}$ of its boundary points and a triangulation
$\mathcal T$  of $\Delta$ whose maximum simplices contain $\bf 0$.
We will combine together the results of Sections
\ref{sec.MP} and \ref{sec.hess} to establish Theorem \ref{main}
which is the main result of this paper. We refer to this theorem as
Toric Residue Mirror Conjecture. We will explain in Section \ref{sec.comp}
that it implies  the original conjecture of \cite{BM}.

To explain the statement of Theorem \ref{main} we need to introduce 
the notion of \emph{toric residues}, as described in \cite{BM}. 
The cone $K$ in the lattice
$\bar M\cong M\oplus \ZZ$ is defined as the span of $\Delta\oplus 1$.
We introduce $\vv_i=v_i\oplus 1$ and $\vv_0={\bf 0}\oplus 1$. 
As before, we consider a generic Laurent polynomial  $f=1+\sum_{i=1}^n a_i t^{v_i}$
and set $a_0=1$. 
Pick a basis $\nn_0,\ldots,\nn_d$ of $\bar N=\bar M^*$. 
The quotient of the graded ring $\CC[K]$ by the elements 
$$
Z_j=\sum_{i=0}^n (\nn_j\cdot \vv_i) a_i t^{\vv_i},~i=0,\ldots, d
$$
is a graded Gorenstein Artin ring, and its degree $d$ component is one-dimensional.
It is spanned by the Hessian $H'_f$ considered in Section \ref{sec.hess}.
\emph{Toric residue} is a map 
$${\rm Res}_f:\CC[K]_{d}\to \CC$$
uniquely defined by its vanishing on the degree $d$ component of the ideal
$\langle Z_0,\ldots, Z_d\rangle \CC[K]$ and by the normalization
$$
{\rm Res}_f(H'_f) = {\rm Vol}(\Delta).
$$
For a given point $\pp=p\oplus d$ in $K$ the value of
${\rm Res}_f(t^\pp)$ is  a rational function in $a_i$ with denominator equal
to the \emph{principal determinant} 
$E=E(a_1,\ldots,a_n)$, see \cite[Theorem 2.9]{BM}.

The principal determinant $E$ is a Laurent polynomial in $a_i$. We can think
of its monomials as being indexed by a lattice $\ZZ^n$ with basis $\{e_i\}$. Vertices 
of the Newton polytope of $E$ (also called \emph{secondary polytope}) are exactly  \emph{characteristic functions} 
$$\chi_{\mathcal T}=\sum_i \sum_{\sigma:v_i\in\sigma\in {\mathcal T},{\rm dim}
\sigma =d} V(\sigma) e_i
$$
which correspond to \emph{regular} triangulations $\mathcal T$ of $\Delta$
whose set of vertices is a subset of $\{{\bf 0}\}\cup\{v_i,i=1,\ldots,n\}$.
The proofs of these statements are contained in \cite{GKZ}, see also \cite{BM}.
Here the triangulation is called regular if there exists a convex
piece-wise linear function $\Delta\to \RR$ whose domains of linearity are precisely
the simplices of $\mathcal T$.

From now on we will assume that the triangulation
${\mathcal T}$ is regular, which also means that the toric variety $\PP_\Sigma$ 
is projective. The rational function ${\rm Res}_f(t^\pp)$ can be expanded
in a Laurent series expansion in the normal cone of the vertex $\chi_{\mathcal T}$ 
of the Newton polytope of $E$, see \cite[Definition 4.5]{BM}.
This cone can be described as follows.
Let $h:M\to \ZZ$ be a convex
function which is linear on the cones of $\Sigma$ and corresponds to an ample
divisor on $\PP_\Sigma$. We extend $h$ to $M_\RR$ by linearity. Then 
for any set of points $y_i \in M_\RR$ and any positive numbers $\alpha_i$
there holds 
\begin{equation}\label{h}
h(\sum_i \alpha_i y_i) \leq \sum_i \alpha_i h(y_i)
\end{equation}
with the equality achieved if and only if there is a cone $\sigma\in \Sigma$ that
contains all points $y_i$. The set of such convex functions will be called the
ample cone of $\mathcal T$ and will be 
denoted by $C^{ample}_{\mathcal T}$.

The following proposition is well-known, but we were unable to find a good
reference in the literature.
\begin{proposition}(see also \cite[Remark 4.7]{BM})
The normal cone $C_{\mathcal T}$ to the vertex $\chi_{\mathcal T}$ in $\ZZ^n$ 
can be characterized by the condition 
$$
\sum_{i=1}^n \alpha_i e_i\in C_{\mathcal T}
 \Leftrightarrow  \sum_{i=1}^n\alpha_i v_i ={\bf 0}
{\rm ~and~}\sum_{i=1}^n \alpha_i h(v_i)\geq 0{\rm ~for~all~}
h\in C^{ample}_{\mathcal T}.
$$
\end{proposition}

\begin{proof}
The normal cone is generated by the differences 
$\chi_{{\mathcal T}_1}-\chi_{{\mathcal T}}$ of the characteristic functions
over all regular triangulations ${\mathcal T}_1$ with the same set 
of vertices.
For every $\mathcal T$ if $\chi_{\mathcal T} =\sum_{i=1}^n \alpha_i e_i$,
then 
$$
\sum_{i=1}^n \alpha_i v_i = \sum_{\sigma\in{\mathcal T},{\rm dim}\sigma=d}
V(\sigma)\sum_{v_i\in \sigma} v_i
$$
is up to a constant 
the baricenter of $\Delta$ and is therefore independent of the
triangulation. Consequently, all differences between 
various $\chi_{\mathcal T}$ satisfy $\sum_{i=1}^n \alpha_i v_i ={\bf 0}$.

For every triangulation $\mathcal T$ and every collection of values 
$(h_0,\ldots, h_n)\in \RR^{n+1}$ 
there is a unique piecewise linear function $h_{\mathcal T}$ on
$\Delta$ which takes values $h_i$ 
on the vertices $v_i$ and $v_0={\bf 0}$ of $\mathcal T$
and is linear on simplices of $\mathcal T$. Moreover, for 
a general collection $(h_i)$ this function $h_{\mathcal T}$ is convex for exactly one
configuration which corresponds to the "bottom" of the convex hull
of $\{v_i\oplus h_i\}\in M_\RR\oplus \RR$. For such $\mathcal T$ and 
$h_{\mathcal T}$ the value $h_{\mathcal T}(p)$ is the smallest  among all
possible values of $h_{{\mathcal T}^1}$ for all triangulations ${\mathcal T}^1$.
The value of  $\sum_{i}\alpha_i h(v_i)$ for a characteristic function of $\mathcal T$ 
is easily seen to equal the integral of $h_{\mathcal T}$ over $\Delta$. Consequently,
it is the smallest of these values among all vertices of the Newton polytope of $E$
if and only if $h$ is convex on $\mathcal T$, and vice versa.
\end{proof}

We are now ready to state our main result.
\begin{theorem}\label{main}
Let $\mathcal T$ be a regular triangulation of a reflexive polytope $\Delta$.
For every $\pp=p\oplus d\in K$ which lies in the lattice spanned 
by $\vv_0,\ldots,\vv_n$ the formal Laurent  series 
$\int_{\mathcal A} \Psi_\pp$ is the expansion of the 
rational function ${\rm Res}_f(t^\pp)$ at the vertex $\chi_{\mathcal T}$ of 
the Newton polytope of $E$.
\end{theorem}

\begin{remark}
Theorem \ref{main} implies that the series 
$\int_{\mathcal A}\Psi_\pp$ is in fact convergent in some open set of 
$(a_1,\ldots,a_n)\in \CC^n$.
\end{remark}

\begin{proof}
First of all, we can restrict our attention to the sublattice of $\bar M$ spanned
by $\vv_i$. All statements about the residues remain unchanged, except
for a possible change in the normalization by the index of the sublattice.
\emph{This is quite different from only looking at the subring of $\CC[K]$ 
generated by $t^{\vv_i}$, the latter may fail to have a one-dimensional
degree $d$ component of the quotient by the ideal $\langle Z_1,
\ldots, Z_d\rangle$.}

Let $\pp=p\oplus d$ be a point in $K$.
\begin{lemma}\label{eff}
Let $\beta\in \ZZ^{n+1}$ satisfy $\sum_{i=0}^n b_i\vv_i = -\pp$ and
$b_0\leq 0$.
Then $\int_{\mathcal A}\Phi_\beta$ is zero unless 
$$
\sum_{i=1}^n b_i h(v_i) + h(p)\geq 0.
$$
for any $h\in C^{ample}_{\mathcal T}$.
\end{lemma}

\begin{proof}
We have $\sum_{i=1}^n b_i v_i + p=0$. We can rewrite it as
$$
\sum_{i,b_i\geq 0} b_i v_i + p = \sum_{i,b_i<0} (-b_i) v_i.
$$
If $\int_{\mathcal A}\Phi_\beta\neq 0$, then all $v_i$ for which
$b_i<0$ lie in a cone of $\Sigma$. We then use \eqref{h} to show that
$$
\sum_{i=1}^n b_i h(v_i) + h(p)= \big(\sum_{i,b_i\geq 0} b_i h(v_i) + h(p) \big)-
\big(  \sum_{i,b_i<0} (-b_i) h(v_i)\big)
$$
$$= 
\big(\sum_{i,b_i\geq 0} b_i h(v_i) + h(p) \big)-h(\sum_{i,b_i\geq 0} b_i v_i + p)
\geq 0.
$$
\end{proof}

\emph{Proof of Theorem \ref{main} continues.} The above lemma implies 
that the formal Laurent series $\int_{\mathcal A}
\Psi_\pp$ are in fact supported in a finite
number of affine shifts of the cone $C_\mathcal T$. The same is true
for the Laurent expansions of ${\rm Res}_f(t^\pp)$. We denote by 
$F_\pp=\int_{\mathcal A}\Psi_\pp-{\rm Res}_f(t^\pp)$ 
the differences and observe that Proposition \ref{rels},
Theorem \ref{hess} and the definition of the toric residue imply that 
\begin{itemize}
\item
For all $\pp_1=p_1\oplus (d-1)$ there holds
$$
\sum_{i=0}^n a_i(\nn\cdot \vv_i) F_{\pp_1+\vv_i} = 0.
$$
\item
$$
\sum_{J\subseteq \{0,\ldots,n\},
\vert  J\vert =d+1} V(J)^2
\,
\big(\prod_{i\in J} a_i\big)
\,
F_{\sum_{i\in J}\vv_i-\vv_0}  = 0.
$$
\end{itemize}

Since for generic $\{a_i\}$ the elements 
$\sum_{i=0}^n a_i(\nn\cdot \vv_i) t^{\pp_1+\vv_i}$ and $H'_f$ generate
$\CC[K]_{{\rm deg}=d}$, the element $t^\pp$ can be written as their linear
combination with coefficients being rational functions in $\{a_i\}$. Consequently,
there is a polynomial $G(a_1,\ldots,a_n)$ such that 
$G(a_1,\ldots,a_n)F_\pp=0$. We remark that a multiplication of a formal
Laurent series by a polynomial is well-defined. We now use the fact that $F$ is
supported in a finite number of affine shifts of $C_{\mathcal T}$. Let $\phi:\ZZ^n
\to \RR$ be a generic linear function which is positive on $C_{\mathcal T}-\{{\bf 0}\}$.
If $F_\pp\neq 0$ then 
there is a term  $c_\alpha a^\alpha$ of $F_\pp$  which has the smallest value
of $h(\alpha)$ among the terms with $c_\alpha\neq 0$. The same can be said
about $G$, and it is easy to see that the product of these terms in $F_\pp G$ 
does not cancel.

We have thus shown that $F_\pp=0$ for all $\pp=p\oplus d$, which proves 
the theorem.
\end{proof}

\begin{corollary}\label{tobm}
Let $P(x_1,\ldots,x_n)\in\QQ[x_1,\ldots, x_n]$ be a degree $d$ polynomial.
Then the Laurent expansion of the toric residue 
$R_P(a)={\rm Res}_f P(a_1t^{\vv_1},\ldots,a_n t^{\vv_n})$
at the vertex $\chi_{\mathcal T}$ is equal to
$$
\sum_{\beta:\sum_{i}b_i\vv_i=0, b_0\leq 0} 
\int_{A_k}P(D_1,\ldots,D_n)D_0^{-b_0} \prod_{i=1}^n D_i^{k-b_i}
\prod_{i=1}^n a_i^{b_i} 
$$
where $k$ is to be taken sufficiently big for each $\beta$.
\end{corollary}

\begin{proof}
The statement follows from Proposition \ref{module}, Theorem \ref{main}
and definitions of $\int_{\mathcal A}\Psi$.
\end{proof}

\section{Connection to the original version of TRMC}\label{sec.comp}
In this section we will establish the connection between Theorem 
\ref{main} and Toric Mirror Symmetry Conjecture of \cite{BM}. 
We will
also remark on the complete intersection case. 

We are working in the notations of the previous section.
\begin{proposition}
Assume that $v_i$ generate $M$. 
Then Corollary \ref{tobm} implies Conjecture 4.6 of \cite{BM}.
\end{proposition}

\begin{proof}
First we observe that 
$\sum_{i=0}^n b_i\vv_i={\bf 0} \in \bar M$ is equivalent to 
$\sum_{i=1}^n b_i v_i ={\bf 0}\in M$ and $-b_0=b_1+\ldots+b_n$.
The argument of Lemma \ref{eff} shows that the summation
in Corollary \ref{tobm} can be taken over the effective classes $\beta$, 
see \cite{BM}.  We also observe that the change in sign of $a_i$
in our notations accounts for the factor $(-1)^d$ of \cite[Conjecture 4.6]{BM}
and introduces an extra factor $(-1)^{b_1+\ldots+b_n}=(-1)^{-b_0}$ to each term.

As a result, it remains to connect
$$
\int_{A_k}P(D_1,\ldots,D_n)(-D_0)^{-b_0} \prod_{i=1}^n D_i^{k-b_i}
$$
with 
$$
\int_{\PP_\beta}P(D_1,\ldots,D_n) (D_1+\ldots+D_n)^{b_1+\ldots+b_n}
\prod_{b_j<0} D_j^{-b_j-1}
$$
where $\PP_\beta$ is defined in \cite[Proposition 3.2 and Definition 3.3]{BM}.

Without loss of generality we can assume that $b_1,\ldots,b_{n-r}\geq 0$
and 
$b_{n-r+1}, \ldots,b_n<0$. We consider the variety 
$\PP_{\Sigma_{(b_1,\ldots,b_{n-r},0,\ldots,0)}}$ of Remark \ref{nondiag}. 
It is given by a fan $\Sigma_\beta$ in the lattice 
$$
M\oplus\bigoplus_{i=1}^{n-r}\ZZ^{b_i}
$$
The vertices of the fan are $v_{i,j}=v_i\oplus e_{i,j}$ for $i\leq n-r, b_i>0$ and
just $v_i$ for $i>n-r$ or $i\leq n-r, b_i=0$.
The cones are given by the condition that the indices $i$ for which all
$v_{i,j}$ are used lie in a cone of $\Sigma$.

We claim that 
\emph{up to a finite index change of  lattice}
 $\PP_\beta$ is isomorphic 
to the toric subvariety  $\PP'_\beta$ in
$\PP_{\Sigma_{(b_1,\ldots,b_{n-r},0,\ldots,0)}}$
which corresponds to the cone generated by $v_i,i>n-r$.
This variety is empty if 
$v_{n-r+1},\ldots,v_n$ do not form a cone in $\Sigma$ and is otherwise
given by the image of the link of $\sigma = {\rm Span}(v_{n-r+1},\ldots,v_n)$
in $\Sigma_{(b_1,\ldots,b_{n-r},0,\ldots,0)}$ modulo the lattice $M_1=
(\QQ v_{n-r+1} +\ldots+ \QQ v_n)\cap M$. The lattice of $\PP'_\beta$ is the quotient
of the lattice $M\oplus\bigoplus_{i=1}^{n-r}\ZZ^{b_i}$ by $M_1$.

We need to recall the definition of $\PP_\beta$ from \cite{BM}.
Consider the lattice 
$$\ZZ(\beta):=\oplus_{i=1}^{n-r}\ZZ_i(\beta)\cong\oplus_{i=1}^{n-r}\ZZ^{b_i+1}$$
with the basis $w_j^{(i)}$. Our notations differ by a switch of $i$ and $j$ from
that of \cite{BM}. Consider the sublattice in 
$\ZZ(\beta)$
defined by the condition
$$
\sum_{i=1}^{n-r} c_i \Big(\sum_{i=0}^{b_i} w_j^{(i)}\Big) = 0
$$
for every solution of $\sum_{i=1}^n c_i v_i = {\bf 0}$.
Then the fan of $\PP_\beta$ lives in the dual $L$ of this sublattice,
which can be thought of as the quotient of the lattice $\ZZ(\beta)^*$ by 
elements 
\begin{equation}\label{gross}
\sum_{i=1}^{n-r} c_i y_i
\end{equation}
for $\sum_{i=1}^n c_i v_i = {\bf 0}$ where $y_i =\sum_{i=0}^{b_i}w_j^{(i)*}$.
Notice that the images of the elements $\frac 1{b_i+1}v_{i,j}, i\leq n-r$ in 
$M_\QQ\oplus\bigoplus_{i=1}^{n-r}\QQ^{b_i}/(M_1)_\QQ$ 
satisfy the same relations \eqref{gross} as the images of $w_j^{(i)}$ in $L$. 
Maximum-dimensional cones of the fan of $\PP_\beta$ are described in the proof  
of Proposition 3.2 of \cite{BM} and it is easy to see that they are in one-to-one
correspondence with the cones of the fan of $\PP'_\beta$. We also remark
that the projection in that proof shows that if $v_{n-r+1},\ldots,v_n$ do not lie in 
a cone of $\Sigma$, $\PP_\beta$ is empty.

Because varieties $\PP_\beta$ and $\PP'_\beta$ are isomorphic
up to a lattice change, their cohomology rings are isomorphic
with the isomorphism mapping 
$D_{i,j}$ to $\frac1{b_i+1}D_j^{(i)}$ which corresponds to 
the element $w_j^{(i)^*}$. As a result, their Stanley-Reisner descriptions in terms of 
the polynomial ring $\CC[D_1,\ldots, D_n]$ have exactly the same ideal,
which includes $D_{n-r+1},\ldots,D_{n}$.
Now we only need to make sure that the top class evaluations are the same, 
which again amounts to an index calculation for some maximum-dimensional cone.
The details are left to the reader.
\end{proof}

\begin{remark}
In general, it appears that Conjecture
4.6 of \cite{BM} needs to be adjusted by the index of the sublattice of $M$ 
generated by $v_i$ inside the lattice $M$. 
For instance, $\PP_{{\bf 0}}$ is in general not isomorphic to $\PP_\Sigma$ but 
is rather a non-ramified abelian cover of $\PP_\Sigma$.
\end{remark}

\begin{remark}
While higher Stanley-Reisner rings $A_k$ are easier to define and work with,
they lack the direct geometric motivation of the toric moduli spaces 
$\PP_\beta$ of \cite{BM}. It is also quite possible that they are better thought
of as Deligne-Mumford stacks, see \cite{BCS}.
\end{remark}

\begin{remark}
It is reasonable to expect that the techniques of this paper are applicable to the complete intersection case of the conjecture, see \cite{BM2}.
\end{remark}

\begin{remark}
It would be interesting to try to apply higher Stanley-Reisner rings to other
open problems in the area. For example one can try to use them to
bound the regularity of the subring of $\CC[K]$ generated by $t^{\vv_i}$.
Surprisingly little is known about this toric case of the more general
Eisenbud-Goto conjecture \cite{EG}.
\end{remark}

\end{document}